\documentstyle[11pt]{amsart}

\title{Area minimizers in a K3 surface and holomorphicity }
\author {M. Micallef and J. Wolfson}
\thanks{The second author was partially supported by
NSF grant DMS-0304587.}

\date{\today}

\newtheorem{thm}{Theorem}[section]
\newtheorem{lem}[thm]{Lemma}
\newtheorem{cor}[thm]{Corollary}

\theoremstyle{definition}

\numberwithin{equation}{section}

\renewcommand{\a}{\alpha}
\renewcommand{\b}{\beta}
\renewcommand{\d}{\delta}

\newcommand{\e}{\varepsilon}
\newcommand{\g}{\gamma}

\renewcommand{\l}{\lambda}

\newcommand{\p}{\partial}

\newcommand{\s}{\sigma}
\newcommand{\Sig}{\Sigma}

\renewcommand{\t}{\tau}
\renewcommand{\O}{\Omega}
\renewcommand{\o}{\omega}
\newcommand{\r}{\rho}

\newcommand{\z}{\zeta}

\newcommand\SI{{\Sigma}}
\newcommand\si{{\sigma}}
\newcommand\xip{{\xi^{1,0}}}
\newcommand\xim{{\xi^{0,1}}}
\newcommand\nup{{\nu^{1,0}}}
\newcommand\num{{\nu^{0,1}}}
\newcommand\pnp{{\pi_{\nu}^{1,0}}}
\newcommand\pnm{{\pi_{\nu}^{0,1}}}

\newcommand\pxm{{\pi_{\xi}^{0,1}}}
\newcommand\De{{D_{\eta}}}
\newcommand\hk{hyperk\"{a}hler}
\newcommand\rd{{\partial}}
\newcommand\cN{{\cal N}}
\newcommand\cU{{\cal U}}
\newcommand\cV{{\cal V}}

\def\RP{\ifmmode{\Bbb R \Bbb P}\else{$\Bbb R \Bbb P$}\fi}
\def\CP{\ifmmode{\Bbb C \Bbb P}\else{$\Bbb C \Bbb P$}\fi}
\def\Z{\ifmmode{\Bbb Z}\else{$\Bbb Z$}\fi}
\def\Q{\ifmmode{\Bbb Q}\else{$\Bbb Q$}\fi}
\def\C{\ifmmode{\Bbb C}\else{$\Bbb C$}\fi}
\def\R{\ifmmode{\Bbb R}\else{$\Bbb R$}\fi}
\def\S{\ifmmode{S^2}\else{$S^2$}\fi}

\def\P{\cal P}

\def\S{\cal S}

\begin{document}

\maketitle

\setcounter{secnumdepth}{1}

\setcounter{section}{0}

\section{\bf Introduction} The following is a  well known consequence of
the Wirtinger inequality: a compact complex submanifold of
a K\"ahler manifold is a volume minimizer in its homology class and
any other volume minimizer in that class is, necessarily, complex.
In particular, in a K\"ahler surface a holomorphic curve is
an area minimizer in its homology class.  In light of this result
it is natural, given a K\"ahler surface, to investigate
the relation between area minimizers and complex curves.
When the ambient manifold is a flat four-dimensional torus
it was shown in [M], using second variation arguments,
that a two-dimensional area minimizer is holomorphic for
one of the complex structures compatible with the metric.
In [MW] the authors attempted to prove an analogous result
for a Ricci flat (Calabi-Yau) metric $g$ on a K3 surface $X$.
Such a metric is hyperk\"ahler in the sense that
there is a two-sphere of complex structures,
called the hyperk\"ahler line, each of which is compatible with $g$.
We obtained partial positive results which are extended
in section 5 of this paper. We also showed that there is
a strictly stable minimal two-sphere in
a (non-compact) hyperk\"ahler surface that is not holomorphic
for any compatible complex structure. This shows that
second variation arguments cannot be used to answer this question.

However there is compelling evidence for a result on K3 surfaces
analogous to that for the flat torus. It can be shown that
if a cohomology class $\a \in H^{1,1}(X ; \R) \cap H^2(X ; \Z)$
satisfies $\a \cdot \a \geq -2$  then its Poincar\'e dual
can be represented by a curve that is holomorphic. Moreover,
there is a set of generators of $H_2(X;\Z)$
each of which can be represented by a curve holomorphic
for some complex structure on the hyperk\"ahler line.
Thus every class $\g \in H_2(X;\Z)$ can be represented by
a sum of curves each of which is holomorphic for
some complex structure on the hyperk\"ahler line.
A minimizer of area  among surfaces representing $\g$
consists of a sum of branched immersed surfaces
$$\Sig_1 \cup \dots \cup \Sig_k,$$
and it is then reasonable to ask whether
each $\Sig_i$ is holomorphic for some complex structure
on the hyperk\"ahler line determined by $g$.
Though this is true for many homology classes we show,
in this paper, that there is an integral  homology class $\g$ and
a hyperk\"ahler metric $g$  such that no area minimizer of $\g$
has this property. Thus the result for flat four-tori
does not carry over to K3 surfaces. For recent work on
similar problems in K\"ahler-Einstein manifolds see [AN].

For lagrangian area minimizers and lagrangian homology classes
there are analogous questions. Given a lagrangian homology class
$\g \in H_2(X ; \Z)$ a minimizer of area
among lagrangian two-spheres representing $\g$ consists of
a sum of lagrangian two-spheres
$$(S^2)_1 \cup \dots \cup (S^2)_k,$$
that may have isolated singularities, as well as branch points.
If each surface is a branched immersion then it can be shown [SW]
that each surface is special lagrangian and, therefore,
holomorphic for some complex structure on the hyperk\"ahler line
determined by $g$. However using techniques similar to
those used in the previous problem it can be shown that
there is an integral lagrangian homology class and
a hyperk\"ahler metric such that no minimizer of area among
lagrangian two-spheres consists solely of branched immersions.
In particular,  there is a lagrangian two-sphere which is
a minimizer of area among lagrangians that is not regular
(i.~e., is not a branched immersion).
We will briefly describe this argument.
A simpler proof of this result, using different techniques,
has been given in [W].
Finally we give a new proof, along the lines of the arguments in [MW],
of a theorem of Donaldson [D] relating stability, holomorphicity and
the normal Euler number of a surface in a K3 surface.
This result suggests that the class we use in
the construction of our main result is optimal in a certain sense.
The work on this paper began at the IPAM workshop:
The Geometry of Lagrangian Submanifolds held at IPAM in April, 2003.
The authors are indebted to IPAM for the hospitality
they extended to us during this workshop.
The first author would also like to thank Mark Gross for
many useful conversations about K3 surfaces.

\bigskip
\section{\bf Preliminaries} In this section we review
basic results in K\"ahler geometry and the geometry of
K3 surfaces that will be used in the proof of our result.
For proofs see [B-P-V] and [G-H].

Let $X$ be a K3 surface, that is, $X$ is a compact, complex,
simply connected surface with trivial canonical bundle. Let
$$L = -E_8 \oplus -E_8 \oplus H \oplus H \oplus H,$$
define the intersection form on a vector space of real dimension 22.
Set $L_\C=L \otimes \C$ with the intersection form extended
complex linearly. For any $\Omega \in L_\C$ we denote
$[\Omega] \in \P(L_\C)$ the corresponding line.
It is known that $H^2(X, \Z)$ is free of rank 22 and
the intersection form on $H^2(X, \Z)$ is given by $L$.
In particular, $b^2_+ = 3$ and $b^2_-=19$.
A marking of $X$ is a choice of basis,
$$ \{ \a_0, \dots, \a_7, \b_0, \dots, \b_7,
\xi_1, \xi_2, \xi_3, \eta_1, \eta_2, \eta_3 \}$$
of $H^2(X, \Z)$ that induces the intersection form $L$.
Equivalently a marking of $X$ is the choice of an isometry
$\phi \colon H^2(X, \Z) \to L$. The period domain ${\cal D}$ of $X$
is the projectivization of the set:
$$\{ \Omega \in L_ \C: \Omega \cdot \Omega =0,\;\;
\Omega  \cdot  \overline{ \Omega} > 0 \}.$$
The complex dimension of ${\cal D}$ equals 20.
If $\Omega$ is a holomorphic $(2,0)$-form on $X$
then the identities $\Omega \cdot \Omega =0$ and
$\Omega  \cdot  \overline{ \Omega} > 0$ show that
a marking of $X$ determines a point $[\Omega] \in {\cal D}$,
called the {\it period point of $X$}.  The first main theorem
we require is the weak Torelli theorem:

\begin{thm}
\label{thm:torelli}
Two K3 surfaces are isomorphic (as complex surfaces) if and only if
there are markings for them such that the corresponding period points
are the same.
\end{thm}

The second main theorem we require is:
\begin{thm}
\label{thm:surj}
All points of the period domain ${\cal D}$ occur as
period points of marked K3 surfaces.
\end{thm}

\medskip
A class $\o \in H^{1,1}(X, \R)$ that can be represented by
a K\"ahler form is called a {\it K\"ahler class}.
Clearly a K\"ahler class satisfies $\o \cdot \o > 0$
and $\o \cdot \O = 0$. Note that the set
$\{ x \in H^{1,1}(X, \R): x \cdot x > 0\}$
consists of two disjoint connected cones and that the K\"ahler classes,
if they exist, all belong to one of these two cones.
This cone is called the {\it positive cone}.
Additional conditions on the K\"ahler classes arise from
the Picard lattice. Let $j \colon H^2(X, \Z) \to H^2(X,\R)$
and define the Picard lattice
$S_X = H^{1,1}(X,\R) \cap \mbox{Im} j(H^2(X,\Z))$.
An element $\s \in S_X$ is called {\it divisorial} if
there exists a divisor $D$ whose associated line bundle
has Chern class $\s$. Then $\s$ is called effective if,
in addition, $D$ can be chosen effective.
The K\"ahler cone is defined to be the convex subcone of
the positive cone consisting of those classes that have
positive inner product with any effective class in $S_X$.
The K\"ahler cone contains all K\"ahler classes.
When $X$ is a K3 surface the characterization of
the K\"ahler cone becomes particularly simple.
A nonsingular curve $\g$ in $X$ is called nodal if
$\g \cdot \g = -2$.

\begin{thm}
For a K3 surface the K\"ahler cone consists of the classes
$\o \in H^{1,1}(X, \R)$ that satisfy: (i)  $\o \cdot \o > 0$,
(ii) $\o \cdot \O = 0$ and (iii) $\o \cdot \g > 0$,
for all nodal curves $\g$ in $X$.
\end{thm}

It is a consequence of the surjectivity of the refined period map
that every class in the K\"ahler cone is a K\"ahler class.
Consequently, Yau's theorem on the existence of
K\"ahler Ricci flat metrics on K3 surfaces can be stated as:

\begin{thm}
Let $(X, \o)$ be a K3 surface where $\o \in H^{1,1}(X, \R)$
lies in the K\"ahler cone.  Then there is
a unique hyperk\"ahler metric on $X$ whose  K\"ahler form represents
the class $\o$.
\end{thm}

\medskip
If $X$ is a K\"ahler surface and $\Sigma$ is a possibly singular
holomorphic curve of genus $g$ in $X$ the {\it adjunction formula} is:
$$\Sigma \cdot \Sigma \geq c_1(X) \cdot \Sigma + 2g - 2,$$
with equality when $\Sigma$ is nonsingular.
When $X$ is a K3 surface this becomes:
$$\Sigma \cdot \Sigma \geq 2g - 2 \geq -2.$$
If $X$ is a K3 surface, we say a (singular) holomorphic curve $\Sigma$
is a $(-2)$-curve if $\Sigma \cdot \Sigma = -2$ (equivalently,
if its Poincar\'e dual $\a$ satisfies $\a \cdot \a = -2$).
From the adjunction formula it follows  that if $\Sigma$ is
a $(-2)$-curve then $\Sigma$ is a nonsingular rational curve.

\bigskip
We conclude this section with some results  on
lagrangian stationary surfaces in K\"ahler-Einstein surfaces
(see [S-W 1]). Let $N$ be a K\"ahler-Einstein surface and
$\Sig$ be a lagrangian submanifold. We say $\Sigma$ is
{\it lagrangian stationary} if the volume is stationary for
arbitrary smooth variations preserving the lagrangian constraint.

\begin{thm}
A closed, branched immersed, lagrangian surface in a
K\"ahler-Einstein surface is a classical minimal surface
if and only if it is lagrangian stationary.
\end{thm}

Consequently,

\begin{cor}
A closed, branched immersed, lagrangian stationary surface in a
K3 surface, with a hyperk\"ahler metric $g$, is  special lagrangian.
In particular, a closed, branched immersed, lagrangian stationary
surface $\Sig$ is a $J$-holomorphic curve with respect to
a complex structure $J$ on the hyperk\"ahler line of $g$.
\end{cor}

Note that these results require regularity of
the lagrangian stationary submanifold.

\bigskip

\section{\bf The Results} In one of the $E_8$'s in
the intersection form of the K3 surface label the four classes
$\a_0, \a_1, \a_2, \a_3 \in H_2(X, \Z)$ that satisfy:
$\a_0 \cdot \a_i = 1$, $i=1,2,3$, $\a_i \cdot \a_j = 0$
for $i \neq j$, $i,j=1,2,3$, $\a_0 \cdot \a_0 = -2$,
and $\a_i \cdot \a_i = -2$.

\begin{lem}
\label{lem:cxstructure}
There is a  complex structure $[\O]$ on the marked K3 surface $X$,
determined by the complex 2-form $\O$, that satisfies:
\begin{enumerate}
\item span$_{\Z} \{\a_0, \a_1, \a_2, \a_3\}
\subset H^{1,1}(X, \C) \cap H^2(X, \Z)$.
\item For each $t$, $0 < t < 1$ there is a K\"ahler class $\o_t$
in the K\"ahler cone determined by $\O$ such that
$\o_t \cdot \a_0 = t$, $\o_t \cdot \a_i = 1$ for $i=1,2,3$
and otherwise $\o_t$ is fixed in $t$.
\item For all $t \in (0,1)$, there exists $\l_t > 0$ such that
$\o_t \cdot \o_t = \frac12  \O_t \cdot \overline{\O_t}$, where
$\O_t := \l_t \O$. Note that $[\O_t] = [\O]$, i.~e.,
the complex structure determined by $\O_t$
is the same as that determined by $\O$.
\item $\sqrt{ (\o_t \cdot \g)^2 + |\O_t \cdot \g|^2 } \geq 1$,
for all $\g \in H^2(X, \Z)$ such that $\g \cdot \g \geq -2$
except when $\g = \pm \a_0$. Equality holds if and only if
$\g = \pm \a_i$ for $i=1,2,3$.
\end{enumerate}
\end{lem}

\begin{pf}
Let a marking of $X$ be given by,
$$ \{ \a_0, \dots, \a_7, \b_0, \dots, \b_7,
\xi_1, \xi_2, \xi_3, \eta_1, \eta_2, \eta_3 \},$$
where the classes $\a_0, \a_1, \a_2, \a_3$ are as described above.
Recall that $\xi_i \cdot \xi_j = 0$, $\eta_i \cdot \eta_j =0$,
for all $i,j$, that $\xi_i \cdot \eta_j = 0$ for $i \neq j$
and that $\xi_i \cdot \eta_i = 1$ for all $i$. Therefore,
$(\xi_i - \eta_i) \cdot (\xi_j + \eta_j) = 0$ for all $i,j$ and
$(\xi_i - \eta_i)^2 = -2$,  $(\xi_i + \eta_i)^2 = 2$ for all $i$.

Define the period point $\O$ as follows:
$$ \O \cdot \a_k = \O \cdot \b_k = 0, \;\;
\mbox{ for } \; k =0, \dots, 7,$$
$$ \O \cdot (\xi_j - \eta_j) = 0, \;\;
\mbox{ for } \; j =1,2,3,$$
$$ \O \cdot (\xi_j + \eta_j) = \s_j + i \t_j, \;\;
\mbox{ for } \; j =1,2,3,$$
where the vectors
$$ \s = (\s_1, \s_2, \s_3), \;\; \t = (\t_1, \t_2, \t_3)$$
satisfy $|\s| = |\t| > 0$ and $\s \cdot \t = 0$. Choose $\s$ so that
no rational linear combination of its components vanishes.
It follows that $\O$ satisfies $\O \cdot \O = 0$ and
$\O \cdot \bar{\O} > 0$ and therefore, by the Torelli theorem,
$\O$ defines a complex structure. By the choice of $\s$
no integral homology class containing a multiple of $\xi_j + \eta_j$
for $j=1,2,3$ can be represented by a holomorphic curve.

Define the K\"ahler class $\o_t$ for $0 < t < 1$ as follows:
$$ \o_t \cdot \a_0 = t, \; \o_t \cdot \a_j = 1, \;
\mbox{ for } \; j =1, 2,3, $$
$$\o_t \cdot \a_j = 2,  \;\;
\mbox{ for} \; j =4, \dots,7,$$
$$ \o_t \cdot \b_j = 2,  \;\;
\mbox{ for } \; j =0, \dots,7, $$
$$ \o_t \cdot (\xi_j - \eta_j) = 2,\;\;
\mbox{ for} \; j =1,2,3,$$
$$ \o_t \cdot (\xi_j + \eta_j) = \r_j,\;
\mbox{ for } \; j =1,2,3,$$
where the vector
$$ \r = (\r_1, \r_2, \r_3) $$
satisfies $\r \cdot \s= \r \cdot \t = 0$. By choosing
$|\r|$ sufficiently large we can ensure that $\o_t \cdot \o_t > 0$.
The classes of nodal curves are integral linear combinations of
$ \{ \a_0, \dots, \a_7, \b_0, \dots, \b_7,
(\xi_1-\eta_1), (\xi_2-\eta_2),( \xi_3-\eta_3) \}$
with coefficients that are non-negative. Therefore
$\o_t \cdot \g > 0$ for all nodal curves. It follows that
$\o_t$ lies in the K\"ahler cone determined by $\O$.

Item (3) is established simply by multiplying $\s$ and $\t$
by an appropriate choice of $\l_t >0$.

For the fourth item, suppose that
$\g \cdot \g \geq -2$ and
$\sqrt{ (\o_t \cdot \g)^2 + |\O_t \cdot \g|^2 } \leq 1$.
We will show that, by taking $|\r|$
(and therefore $|\s_t|$ and $|\t_t|$) sufficiently large,
$\g$ then has to be one of $\pm \a_0,\ \pm \a_1,\
\pm \a_2,\ \pm \a_3$.

Decompose $\g$ into its self-dual and anti self-dual parts
$\g_+$ and $\g_-$. Thus $\g_+ = \sum_{i=1}^3 n_i(\xi_i+\eta_i), \
n_i \in \frac12\Z$ and $\g_-\in\text{span}_{\Z}
\{\a_0, \dots, \a_7, \b_0, \dots, \b_7,
\frac12(\xi_1-\eta_1), \frac12(\xi_2-\eta_2), \frac12(\xi_3-\eta_3) \}$.
As a first step, we shall show that $\g_+$ has to vanish
when $|\r|$ is sufficiently large.
Let $\langle\cdot,\cdot\rangle$ and $|\cdot|$ denote
the Euclidean inner product and norm on $\R^3$. Then,
letting $n := (n_1,n_2,n_3)$ we get:
\begin{equation}\label{g2}
-2 \leq \g\cdot\g = 2|n|^2 + (\g_-\cdot\g_-),
\end{equation}
where
\begin{equation}\label{n2}
|n|^2 = \frac{|\langle n,\s_t \rangle|^2}{|\s_t|^2} +
\frac{|\langle n,\t_t \rangle|^2}{|\t_t|^2} +
\frac{|\langle n,\r \rangle|^2}{|\r|^2}.
\end{equation}
We need to show that $n=(0,0,0)$. If not, then $|n|^2 \geq 1/4$.
Furthermore, by assumption, $|\O_t \cdot \g| \leq 1$ and so,
$|\langle n,\s_t \rangle| \leq 1/2$ and
$|\langle n,\t_t \rangle| \leq 1/2$.
Putting these inequalities in (\ref{n2}) yields
$$|\langle n,\r \rangle|^2 \geq |n|^2|\r|^2
\left(1-\frac{2}{|\s_t|^2} - \frac{2}{|\t_t|^2}\right). $$
It will be useful to observe that
$|\s_t|^2 = |\t_t|^2 = \o_t \cdot \o_t$ for all $t \in (0,1)$
and that there exists $C>0$ such that for all $t \in (0,1)$,
\begin{equation}\label{otr}
\tfrac12 |\r|^2 \geq \o_t \cdot \o_t \geq \tfrac12 |\r|^2 - C.
\end{equation}
It follows that, for sufficiently large $|\r|$,
\begin{equation}\label{nr}
|\langle n,\r \rangle|^2 \geq \tfrac12|n|^2|\r|^2.
\end{equation}
Now $\o_t \cdot \g = 2\langle n,\r \rangle + \o_t \cdot \g_-$.
By our assumption $|\o_t \cdot \g| \leq 1$.  Therefore
$|\o_t \cdot \g_-| \geq 2|\langle n,\r \rangle| - 1$.
From this we deduce
$$ (\o_t \cdot \o_t)(-\g_- \cdot \g_-) \geq
(\o_t \cdot \g_-)^2 \geq (2|\langle n,\r \rangle| - 1)^2
\geq 3(|\langle n,\r \rangle|^2 - 1). $$
Substituting into (\ref{g2}) and making use of (\ref{nr})
and (\ref{otr}) yields:
$$(|n|^2+1)|\r|^2 \geq 2(|n|^2+1)(\o_t \cdot \o_t) \geq
(-\g_-\cdot\g_-)(\o_t \cdot \o_t)
\geq \tfrac32|n|^2|\r|^2 - 3.$$
Clearly this cannot hold for arbitrarily large $|\r|$ unless
$n = (0,0,0)$.

Now that we have shown that $\g_+ = 0$, we see that
$-2 \leq \g \cdot \g = \g_- \cdot \g_- \leq -2$.
It follows that $\g \cdot \g = -2$ and that
$\g \in \text{span}_{\Z} \{\a_0, \dotsc, \a_7\}$ or
$\g \in \text{span}_{\Z} \{\b_0, \dotsc, \b_7\}$.
($\g$ cannot be one of $\pm(\xi_i-\eta_i), \ i=1,2,3$
because $\o_t \cdot (\xi_i-\eta_i) = 2$ by construction
and $|\o_t \cdot \g| \leq 1$ by assumption.)
Suppose that $\g \in \text{span}_{\Z} \{\a_0, \dotsc, \a_7\}$.
A tedious calculation shows that if
$\g = \sum_{i =0}^7 m_i \a_i, \ m_i \in \Z$ and $\g \cdot \g = -2$
then $m_0, \dotsc , m_7$ {\em must all have the same sign}.
Actually, there is a theoretical reason for this.
The set $\{\g \in \text{span}_{\Z} \{\a_0, \dotsc, \a_7\} :
\g \cdot \g = -2\}$ is the root system of type $-E_8$;
there are 240 such roots. (See, for example, [H], pages 472 and 473.)
$\{\a_0, \dotsc, \a_7\}$ is a basis of this root system  and therefore,
if $\g$ is a root, then $\g = \sum_{i =0}^7 m_i \a_i, \
m_i \in \Z$ and $m_0, \dotsc , n_7$ {\em all have the same sign}.
The same applies if $\g \in \text{span}_{\Z} \{\b_0, \dotsc , \b_7\}$.
Item (4) now follows easily.
\end{pf}

\noindent {\bf Remark:} The class $\a_4 - \a_5$ shows that
(4) does not hold if $\g \cdot \g \leq -4$.

\bigskip
By Yau's theorem for each $0 < t < 1$ there is a unique
hyperk\"ahler metric $g_t$ with K\"ahler form in the class of $\o_t$.
Note that the complex structure $[\O]$ is fixed.

Using the Riemann-Roch theorem it follows that
the indecomposable classes $\a_0, \a_1, \a_2, \a_3$
can be represented by embedded holomorphic $-2$ curves
that we denote, respectively $S_0, S_1, S_2, S_3$, [B-P-V, VIII 3.6].
These curves do not change with $t$.

The sequence of Ricci flat K\"ahler metrics $\{g_t \}$ on $X$,
as $t \to 0$, has been studied by R.~Kobayashi [K].
Denote the orbifold obtained from $X$ by blowing down $S_0$ by $X_0$.
Then $X_0$ has one orbifold point $p$ and
an orbifold K\"ahler Ricci flat metric $g_0$
which is singular only at $p$. (The existence of
$g_0$ was established in [KT], Theorem 1, p.348.)
Furthermore, according to Theorem 21 in [K], the metrics
$g_t$ converge smoothly on compact subsets of
$X \setminus S_0$ to the metric $g_0$ on $X_0 \setminus \{p\}$.

Let $\cN$ be a tubular neighborhood of $S_0$ in $X$ and let
$\bar{S_i} := \pi (S_i \cap \cN)$ where
$\pi \colon X \to X_0$ is the blow down projection.
Of course, $\bar{S_1}, \ \bar{S_2}$ and $\bar{S_3}$
all meet at $p$ in $X_0$. In order to understand this intersection
fully, we shall recall the explicit description of $\pi$.
The tubular neighborhood $\cN$ is biholomorphic to a
neighborhood of the zero section of $T^*\CP^1$;
we shall therefore identify this zero section with $S_0$.
The blow down of $T^*\CP^1$ along $S_0$ is the quadratic cone
$\C^2/\Z_2$, where $\Z_2$ acts on $(z_1,z_2)$ by
$(z_1,z_2) \mapsto (-z_1,-z_2)$. Let $\s \colon \C^2 \to \C^2/\Z^2$
be the natural projection. We shall exhibit
a holomorphic double covering
$\r \colon \C^2 \setminus \{(0,0)\} \to T^*\CP^1 \setminus S_0$
such that $\pi \circ \r = \s$. For this purpose, cover $T^*\CP^1$
by two coordinate charts $\cU \cong \cU' \cong \C^2$ with
coordinates $(u,\xi), \ (u',\xi')$ respectively;
here $u$ and $u'$ denote Euclidean coordinates on the base $\CP^1$,
and $\xi$ and $\xi'$ parametrize the fibers of the bundle
$T^*\CP^1 \to \CP^1$. This means that
$$ \cU \cap \cU' = \C^* \times \C \text{ and, if }
(u',\xi') \in \cU \cap \cU' \text{ then }
(u',\xi') \sim (1/u, u^2 \xi). $$
Define $\r \colon \C^2 \setminus \{(0,0)\} \to
T^*\CP^1 \setminus S_0$ by
$$ \r(z_1,z_2) = \begin{cases}
(z_1/z_2,z_2^2) \in \cU \setminus \{\cU \cap S_0 \}, &
\text{if $z_2 \neq 0$,} \\
(z_2/z_1,z_1^2) \in \cU' \setminus \{\cU' \cap S_0\}, &
\text{if $z_1 \neq 0$.} \end{cases} $$
Note that $\r(z_1,z_2) = \r(-z_1,-z_2)$ and therefore,
$\r$ descends to a biholomorphic map $\bar{\r} \colon
\C^2 \setminus \{(0,0)\}/\Z_2 \to T^*\CP^1 \setminus S_0.$
It follows that $\pi(q) = \bar{\r}^{-1}(q)$ if $q \notin S_0$ and
$\pi(q) = [(0,0)]$ if $q \in S_0$; this establishes $T^*\CP^1$
as the minimal resolution of the cone $\C^2/Z_2$.

Let $\hat{S_i} := \s^{-1}(\bar{S_i}) \subset \C^2$.
We shall show that $\hat{S_1}, \ \hat{S_2}$ and $\hat{S_3}$
intersect pairwise transversally at $(0,0)$.
Let $p_i := S_i \cap S_0$. We may as well assume that
$p_i = (u_i,0) \in \cU$. Since $S_i$ intersects $S_0$ transversely
at $p_i$, there exist holomorphic functions $f_i \colon B_{\e} \to \C$,
$B_{\e} := \{w \in \C : |w| < \e\}$, such that
$S_i \cap \cN = \{(f_i(w),w) : w \in B_{\e}\}$. It follows that
$\hat{S_i} = \{(\z f_i(\z^2),\z) : |\z^2| < \e \}$.
But $f_i(0) = u_i$ and $u_1, \ u_2, \ u_3$ are all distinct.
Therefore, $\hat{S_1}, \ \hat{S_2}$ and $\hat{S_3}$ are all nonsingular
and intersect pairwise transversally at $(0,0)$.

To say that $g_0$ is a K\"ahler Ricci flat orbifold metric on $X_0$
means that $\hat{g_0} := \s^*(g_0|_{\cV}), \ \cV := \pi(\cN)$,
is a smooth $\Z_2$ invariant metric on $\s^{-1}(\cV)$, which is a
neighbourhood of $(0,0) \in \C^2$.
Two of $\hat{S_1}, \ \hat{S_2}$ and $\hat{S_3}$
must intersect non-orthogonally with respect to $\hat{g_0}$.
Renumbering we can suppose that this pair is $\hat{S_1}, \ \hat{S_2}$.
Reverse the orientation on $\hat{S_2}$.
Then the two tangent planes of $\hat{S_1}$ and $\hat{S_2}$ at $(0,0)$
intersect non-orthogonally at $(0,0)$ and
define an orientation on $\C^2$ which is opposite to the canonical one.
Hence they do not form an area minimizing configuration in
$T_{(0,0)}(\cV)$ [L], [Mo].
It follows that in $\cV$, there are discs $D_i \subset \hat{S_i}$
centered at $(0,0) \in \hat{S_1} \cap \hat{S_2}$ and
an annulus $A$ in $\cV$ (a ``Lawlor neck'') with boundary
$\p A = \p D_1 \cup \p D_2$ such that,
$$\mbox{area}_{\hat{g_0}} (A) <
\mbox{area}_{\hat{g_0}} (D_1 \cup  D_2).$$
Note that the annulus can be chosen so that
$A$ is disjoint from $(0,0)$. Then for $\e > 0$ sufficiently small
it remains true that,
$$\mbox{area}_{\hat{g_0}} (A) <
\mbox{area}_{\hat{g_0}} ((D_1 \cup  D_2) \setminus
(D_1 \cup  D_2) \cap B_{\e}).$$
Since $\s$ restricted to $\cV \setminus B_{\e}$ is a local isometry,
it follows that in $X_0 \setminus \{p \}$:
\begin{eqnarray*}
\mbox{area}_{g_0} (\s(A)) & < & \mbox{area}_{g_0}
(\s((D_1 \cup  D_2) \setminus (D_1 \cup  D_2) \cap B_{\e}))\\
& < & \mbox{area}_{g_0} (\s(D_1 \cup  D_2)).
\end{eqnarray*}
The annulus $\s(A)$ can then be used to glue
$\bar{S_1} \setminus \s(D_1)$ to $-(\bar{S_2}  \setminus \s(D_2))$
forming a piecewise $C^1$ two-sphere
$\bar{S}  \subset X_0 \setminus \{p\}$.
Clearly $\bar{S}$ has $g_0$-area strictly less than
$\mbox{area}_{g_0}(\bar{S_1}) + \mbox{area}_{g_0}(\bar{S_2}) =2$.
Since the metrics $g_t$ converge in $C^{\infty}$ uniformly
on compact subsets of $X_0 \setminus \{p\}$ to $g_0$,
for $t$ sufficiently small, the two-sphere
$S := \bar{\r}(\bar{S}) \subset X \setminus S_0$
has $g_t$-area strictly less than $2$ in $(X, g_t)$.
In conclusion, there is a two-sphere $S$ in $(X, g_t)$
that represents $\a_1 - \a_2$ in integral homology such that, for
sufficiently small $t$, the $g_t$-area of $S$ is strictly less than $2$.
Recall that the $g_t$-area of each of the holomorphic curves $S_1$,
$S_2$ and $S_3$ is $1$ and of the holomorphic curve $S_0$ is $t$.

\medskip
We say a surface $\Sig$ in the K3 surface $(X, g)$ is
{\it ${\cal J}$-holomorphic} if it is J-holomorphic for some
$J \in {\cal J}$ the hyperk\"ahler line of $g$.
We will need the following elementary lemma.

\begin{lem}
\label{lem:area}
Suppose that $\Sigma$ is a surface in the K3 surface $(X, g_t)$
that is ${\cal J}$-holomorphic. Then
$$\mbox{area}_{g_t} (\Sigma) =
\sqrt { (\o_t \cdot \g)^2 + |\O_t  \cdot \g|^2}$$
where $\g$ is the Poincar\'e dual of $\Sigma$.
\end{lem}

\begin{pf}
Suppose that $\Sigma$ is holomorphic with respect to
the complex structure $J$ compatible with $g_t$ and let $\o$
be the K\"ahler form of $(J, g_t)$. Then
$$\mbox{area}_{g_t} (\Sigma) = \int_{\Sigma} \o = \g \cdot \o.$$
Let $\g_+$ be the self-dual part of $\g$. Then
$(\g \cdot \o)^2 = (\g_+ \cdot \g_+) (\o \cdot \o)$. But
$$\g_+ \cdot \g_+ = \frac{(\g \cdot \o_t)^2}{\o_t \cdot \o_t}
+ \frac{2|\g \cdot \O_t|^2}{\O_t \cdot \overline{\O_t}}.$$
The result follows using $\o_t \cdot \o_t =
\frac12 \O_t \cdot \overline{\O_t}$ and $\o \cdot \o = \o_t \cdot \o_t$.
\end{pf}

The following theorem is our main result.

\begin{thm}
\label{thm:main1}
When $t$ is sufficiently small,
no area minimizer of $\a_1 - \a_2 \in H_2(X, \Z)$,
for the hyperk\"ahler metric $g_t$,
is  the sum of surfaces each of which is ${\cal J}$-holomorphic.
\end{thm}

\begin{pf} Suppose, by way of contradiction, that some area minimizer
of $\a_1 - \a_2$ is the sum of surfaces each of which is
${\cal J}$-holomorphic. Recall that we have constructed a two-sphere
representing this class with area strictly less than $2$.
First suppose that the area minimizer has one component $C$.
Then $C$  is a branched minimally immersed surface that
represents $ \a_1 - \a_2$. Thus,
$$C \cdot C = (\a_1 - \a_2)  \cdot (\a_1 - \a_2) = -4.$$
$C$ cannot be holomorphic, for any complex structure,
as this would contradict the adjunction formula.

Next suppose that the area minimizer has, at least,
two component surfaces $B$ and $C$. Each one is, by assumption,
holomorphic for some complex structure. Therefore,
provided that neither surface represents $\pm \a_0$ in homology,
by Lemmas \ref{lem:cxstructure} and \ref{lem:area},
$$\mbox{area} (B) + \mbox{area} (C) \geq 2.$$
But the sum of the areas of all component surfaces is
stricly less than $2$. Therefore, at least one of the surfaces
represents either $ \a_0$ or $-\a_0$ in homology.
The sum of the remaining surfaces then represents either
(i) $\a_1- \a_0 - \a_2$ or (ii) $\a_1 +  \a_0 - \a_2$.
Note that a single holomorphic curve cannot represent
$\a_1 + k \a_0 - \a_2$ for any $k \in \Z$
(since $(\a_1 + k \a_0 - \a_2)^2 = -4 - 2k^2$). It follows that
in both cases the set of remaining surfaces contains at least
two components neither of which represents $\pm \a_0$ in homology.
The previous argument shows that this is impossible.
The result follows.
\end{pf}

\bigskip
The previous arguments can be adapted to the study of minimizers
among lagrangian two-spheres as follows:
Recall that for each $0 < t < 1$ there is a unique
hyperk\"ahler metric $g_t$ with K\"ahler form in the class of $\o_t$.
Then for each $0 < t < 1$ there is an $S^1$ of K\"ahler forms,
compatible with $g_t$, such that the $-2$-curves $S_1$ and $S_2$
(used above) are lagrangian with respect to each form in the family.
For each $t$ choose such a K\"ahler form and denote it $w_t$.
Choosing a subsequence as $t \to 0$, we can suppose that
the sequence $\{w_{t} \}$ converges to
the (orbifold) symplectic form $w_0$ on $X_0$ and that
$\bar{S_1}$ and $\bar{S_2}$ are lagrangian for $w_0$.
The oriented lagrangian surfaces $\bar{S_1}$ and $-\bar{S_2}$
can  be glued using a Lawlor neck to construct
a lagrangian two-sphere $\bar{S}$ in $X_0 \setminus \{p\}$ ---
$p$ being the orbifold singularity of $X_0$ ---
with $g_0$-area equal to $2 - \e$, for some $\e > 0$.
It follows that there is a $\d > 0$ such that for all $t < \d$
the two-sphere $S := \bar{\r}(\bar{S}) \subset X \setminus S_0$
has $g_t$-area strictly less than $2 - \e/2$ in $(X, g_t)$.
Now on any fixed compact subset of $X \setminus S_0$
$w_{t}$ converges smoothly to $w_0$ . Therefore,
for $\d$ sufficiently small and $t < \d$,
$S$ is approximately lagrangian with respect to $w_t$.
In particular, it can be deformed to a lagrangian two-sphere
without changing its $g_{t}$-area by  more than $\e/2$.
We conclude that for some $t$ sufficiently small
there is a  lagrangian two-sphere with $g_{t}$-area less than $2$.
Next we minimize area among lagrangian two-spheres that represent
$\a_1 - \a_2$.

\begin{thm}
\label{thm:main2}
There is some lagrangian class in $H_2(X; \Z)$ that has
an area minimizer among lagrangian two-sphere
(for the hyperk\"ahler metric $g_t$, $t$ sufficiently small)
that is not a branched immersion.
\end{thm}

\begin{pf}
Suppose, by way of contradiction, that every lagrangian two-sphere
that is a $g_{t}$-area minimizer is regular (i.~e., is
a branched immersion). Then by the results of [SW]
every such lagrangian two-sphere is holomorphic for some
complex structure on the hyperk\"ahler line determined by $g_{t}$.
Consider an area minimizing sequence of lagrangian two-spheres that
represent $ \a_1 - \a_2$ and recall that we have constructed a
lagrangian two-sphere representing this class with area
strictly less than $2$. The argument in the proof of
Theorem \ref{thm:main1} leads to a contradiction.
\end{pf}

\noindent {\bf Remark:} It follows from Theorem \ref{thm:main2} that,
on lagrangian area minimizers in a K3 surface,
singularities other than branch points can and do occur.

To put this result in context, recall the constrained
variational theory developed in [SW]. Consider a homology class
in a K\"ahler surface that can be represented by the image of
a lagrangian map of a compact surface (a {\it lagrangian
homology class}) and minimize area among such maps. Then in [SW]
it is shown that a lagrangian minimizer exists,
that the map is Lipschitz and is an immersion except at
a finite number of isolated points that are either
(i) branch points, or (ii) singular points with non-flat tangent cone.
The tangent cones can be described precisely and it can be shown that
there is a Maslov index associated to each tangent cone
(and hence to each singular point). If the map is a minimizer
this index is $\pm 1$. The sum of these indices equals the pairing of
the first Chern class of the K\"ahler surface with
the homology class of the minimizer. Thus,
when this pairing is non-zero, a lagrangian minimizer must
have singular points. However if the K\"ahler surface is
K\"ahler-Einstein then this pairing vanishes and
it is possible that the minimizer is always regular.
More precisely, one could speculate that on a minimizer
a pair of singularities with indices $1$ and $-1$
could be shown to ``cancel''.
Theorem \ref{thm:main2} shows that this is not the case.

\bigskip
\section{\bf The Motivating Idea} The main construction in this paper
is motivated by the following observation. Suppose a K3 surface $X$
admits two $-2$-curves $S_1$ and $S_2$ that intersect.
Denote $[S_1] = \a_1$ and $[S_2] = \a_2$.
Suppose there is a hyperk\"ahler metric $g$ on $X$ with
K\"ahler form $\o$ such that $\o(S_1) = \o(S_2) = 1$ and such that
$S_1$ and $S_2$ intersect non-orthogonally with respect to $g$.
Then after gluing in a ``Lawlor neck'' a representative of
the class $\a_1 - \a_2$ can be constructed with $g$-area less than $2$.
(In fact, this representative can be taken to be lagrangian
for the K\"ahler form $\mbox{Re}\, \O$.)
Suppose an area minimizing sequence among two-spheres representing
$\a_1 - \a_2$ converges without bubbling. Then, by
the adjunction formula, the area minimizer cannot be holomorphic.
The area minimizing sequence cannot bubble into two spheres representing
$\a_1$ and $-\a_2$, respectively, because then the sum of
the areas of the bubbles is $2$. Other possible bubbling can be
ruled out by appropriate choice of complex structure and metric on $X$,
as was done in Section 3.

We did not use this construction because of the technical difficulty of
finding a hyperk\"ahler metric on $X$ such that two $-2$-curves
intersect non-orthogonally. Our construction in Section 3
exploits the existence of two $-2$-curves that intersect
non-orthogonally in an orbifold limit  of the hyperk\"ahler metrics
on $X$. That is, we find the necessary $-2$-curves on
the boundary of the moduli space of Calabi-Yau metrics
rather than in the interior.

\bigskip
\section{\bf Stability and Holomorphicity}
In this section we give a new proof, based on ideas in [MW],
of a theorem of Donaldson [D]. Donaldson's result implies that
an immersed area minimizer in a K3 surface with normal Euler number
greater than $-4$ must be ${\cal J}$-holomorphic. If the area minimizer
produced in Theorem \ref{thm:main1} is immersed and
consists of one component (i.e., there is no bubbling)
then Donaldson's result shows that this example is optimal.

\begin{thm}\label{don} \emph{(Donaldson)}
Let $\SI$ be an oriented immersed minimal surface
in a 4-manifold $X$ equipped with a \hk\ metric $g$.
If the Euler number $e(\nu)$ of the normal bundle $\nu$
of $\Sig$ in $X$ satisfies $e(\nu) \geq -3$ and
$\Sig$ is not holomorphic with respect to a K\"{a}hler structure
of $g$ then $\Sig$ cannot be strictly stable.
\end{thm}

\begin{pf}
We start by recalling some of the basic geometry of
an immersion $F\colon \Sig \to X$ of an oriented surface $\Sig$
in an oriented 4-manifold $X$ equipped with a Riemannian metric $g$.
For more detail, refer to [MW].
We shall adopt the following notation:
\begin{enumerate}
\item[(i)]$(x,y)$ will denote isothermal co-ordinates for
the metric on $\Sig$ induced by $F$.
$z=x+iy$ is then a local complex co-ordinate on $\Sig$.
\item[(ii)]$F^*(TX) = \xi\oplus\nu, \quad
\xi=\text{tangent bundle of }\SI$, $\nu=\text{normal bundle.}$
\item[(iii)]\hfill $\xi_{\C}:=\xi\otimes_{\R}\C, \quad
\nu_{\C}:=\nu\otimes_{\R}\C,$ \hfill \mbox{}

$\xi_{\C} = \xip\oplus\xim, \quad \xip(\xim) =
+i(-i)$ eigenspace of the rotation $J_{\SI}$ in $\xi$
by $90^{\circ}$ anticlockwise. Similarly,
$\nu_{\C} = \nup\oplus\num$; rotation by $90^{\circ}$ anticlockwise in
$\nu$ is denoted by $J_{\nu}$.
\item[(iv)]Superscript $\top (\perp)$ will denote orthogonal projection
onto $\xi(\nu)$.
\item[(v)]$D\colon\Gamma(F^*(TX) \to \Gamma(T^*(\SI)\otimes F^*(TX)$
is the connection on the pull-back of the tangent bundle of $X$
induced by the Levi-Civita connection of $g$.
We shall also make use of
$D':=dz \otimes D_{\rd/\rd z}$ and
$D'':=d \bar{z} \otimes D_{\rd/\rd \bar{z}}$.
This notation is slightly different from that in
the Appendix in [MW] but it is more conventional.
\item[(vi)]$D|_{\Gamma(\nu)} = \nabla+B;
\quad \nabla:=D^{\perp}|_{\Gamma(\nu)}$
is a metric compatible connection on $\nu$ and
$B$ is the shape operator.
$B'$ is naturally defined by
$B's:=(D's)^{\top}\ \forall\,s\in \Gamma(\nu)$.
$B'',\ \nabla'$ and $\nabla''$ are defined similarly.
\item[(vii)]$D|_{\Gamma(\xi)} = D^{\top}-B^*;\quad D^{\top}$
is the Levi-Civita connection on $\SI$ and
$B^*$ is the second fundamental form of $F$;
it is dual to the shape operator $B$.
Of course, $B^*{}'$ is defined by
$B^*{}'s:=-(D'\si)^{\perp}\ \forall\,\si\in \Gamma(\xi)$.
$B^*{}''$ is defined similarly.
\item[(viii)]$F$ is minimal if, and only if, $B^*{}''|_{\xip} \equiv 0$,
which is equivalent to $\pxm \!\circ\! B'' \equiv 0$.
\end{enumerate}

We now assume that $g$ is \hk\ and we let $I_1,\,I_2,\,I_3$
be complex structures on $X$ which are parallel with respect to
the Levi-Civita connection of $g$ and which define
a quaternionic structure on $TX$ (and therefore, $F^*(TX)$).
Let $J:=J_{\SI}\oplus J_{\nu}.$ Then
$J = \sum_{j=1}^3(\cos \alpha_j)I_j,\
\sum_{j=1}^3\cos^2 \alpha_j \equiv 1$.
$\alpha_j$ is called the K\"ahler angle of $\SI$ with respect to $I_j$.
A straightforward calculation shows that, for each $j \in \{1,2,3\}$,
$I_j$ preserves $\Gamma(\nup\oplus\xim)$.
This is essentially due to the fact that Hom$(\nup,\xim)$
is the tangent space of the twistor space of $F^*(TX)$ at $J$.
Let $\eta:=\nup\oplus\xim$ and let $\pi_{\eta} = \pnp\oplus\pxm$.
Define a connection $\De$ on $\Gamma(\eta)$ by
$\De:=\pi_{\eta}\!\circ\!D|_{\Gamma(\eta)}$.
Endow $\eta$ with a holomorphic structure by declaring
$v \in \Gamma(\eta)$ to be holomorphic if, and only if,
$\De''v \equiv 0$ where
$\De'':=\pi_{\eta} \!\circ\! D''|_{\Gamma(\eta)}$.
$\nup$ is endowed with the holomorphic structure
defined by means of $\nabla''$ and
$\xim$ is endowed with the holomorphic structure
defined by means of $D''{}^{\top}$. Assume that $F$ is minimal and
write $v=s+\si,\ v\in\Gamma(\eta),\ s\in\Gamma(\nup),\
\si\in\Gamma(\xim)$. Then
$\De''v = \nabla''s - B^*{}''\si + D''{}^{\top}\si$
where we have used $\pxm\!\circ\!B'' \equiv 0$
by the minimality of $F$. Thus,
\[
v\in H^0(\eta) \Leftrightarrow \si\in H^0(\xim) \quad
\text{and} \quad \nabla''s = B^*{}''\si
\]
where $s=\pnp v$ and $\si=\pxm v$.
The proof of the theorem makes use of
the second variation of area formula applied to
sections $s=\pnp v,\ v\in H^0(\eta)$. In [MW],
the second variation of area formula was applied to
holomorphic sections of $\nup$, which can also be viewed as
holomorphic sections of $\eta$ (with $\si=0$).

Let $\cal R$ denote the curvature operator of $X$ and
let $e_1-ie_2$ be a local unitary section of $\xip$.
Then, according to (A.10) in [MW],
the second variation $(\delta^2A)(s)$ of area
in the direction of $s\in\Gamma(\nu_{\C})$
of a minimal surface $\SI$ in $X$ is given by:
\[
(\delta^2A)(s) = 2\int_{\SI}\{|\nabla''s|^2-|B's|^2
-\tfrac12\langle{\cal R}s\wedge(e_1-ie_2),s\wedge(e_1-ie_2)\rangle\}
\,dA.
\]
If $s\in\Gamma(\nup)$ and $X$ is K\"{a}hler with zero scalar curvature
then, as can be seen from Proposition 2.2 in [MW],
the curvature term in the formula for $(\delta^2A)(s)$ drops out.
Furthermore,
\begin{equation}\label{Bprime}
|B's|^2 = -\tfrac12(K_{\xi}+K_{\nu})|s|^2
\end{equation}
where $K_{\xi}$ is the Gauss curvature of $\SI$ and
$K_{\nu}$ is the curvature of the connection $\nabla$ on $\nu$.
This can be seen as follows: let $s=f(e_3-ie_4)$ where
$e_3-ie_4$ is a local unitary section of $\nup$. Then
\[
B's = \tfrac12 f(h_{311}-h_{412}-ih_{312}-ih_{411})
(e_1+ie_2)\otimes(\theta_1 + i \theta_2)
\]
where $\theta_1 + i \theta_2$ is the local unitary section of
$(\xip)^*$ dual to $e_1-ie_2$. Therefore,
\[
|B's|^2 = \tfrac12 |s|^2
\bigl((h_{311}-h_{412})^2 + (h_{312}+h_{411})^2 \bigr).
\]
Let $R_{ABCD}:=\langle{\cal R}(e_A \wedge e_B),e_C \wedge e_D \rangle$.
Then, by the Gauss equation and the minimality of $F$,
\[
K_{\xi} = R_{1212} - (h_{311}^2+h_{312}^2+h_{411}^2+h_{412}^2)
\]
and, by the Ricci equation,
\[
K_{\nu} = R_{3412} + 2(h_{311}h_{412}-h_{312}h_{411}).
\]
Equation (\ref{Bprime}) now follows easily on noting that
$R_{3412}+R_{1212} =
\langle{\cal R}(e_1 \wedge e_2 + e_3 \wedge e_4),e_1 \wedge e_2 \rangle$
which, according to Proposition 2.2 in [M-W], is equal to zero
for a K\"{a}hler surface $X$ with zero scalar curvature.

We are now in a position to prove Theorem \ref{don}.
Let $v$ be a nontrivial holomorphic section of $\eta$ and let
$s_j:=\pnp(I_jv),\ j\in\{0,1,2,3\}$ where, for convenience,
$I_0$ denotes the identity transformation on $\Gamma(\eta)$.
Then, since for each $j,\ I_jv$ is a holomorphic section of $\eta$,
we have $\nabla''s_j = B^*{}''\si_j$ where $\si_j:=\pxm(I_jv)$.
Now by a calculation similar to that just carried out
for the establishment of equation (\ref{Bprime}) one can show that
\[
|B^*{}''\si_j|^2 = -\tfrac12(K_{\xi}+K_{\nu})|\si_j|^2.
\]
It follows that
\[
\sum_{j=0}^3(\delta^2A)(s_j) =
\int_{\SI}(K_{\xi}+K_{\nu})\sum_{j=0}^3(|s_j|^2-|\si_j|^2)\,dA\,.
\]
It is easy to see that
\[
\sum_{j=0}^3|s_j|^2 = \sum_{j=0}^3|\si_j|^2 = 2|v|^2.
\]
Therefore, $\sum_{j=0}^3(\delta^2A)(s_j) = 0$.
But there are at least two values of $j$
(which we may as well assume to be 0 and 1)
for which $s_j$ does not vanish identically.
Strict stability would then imply $(\delta^2A)(s_j) > 0,\
j \in \{0,1\}$ and $(\delta^2A)(s_k) \geq 0,\ k \in \{2,3\}$.
Hence strict stability and the existence of a nontrivial
holomorphic section of $\eta$ lead to a contradiction.

The proof of the theorem is completed by showing that,
if $e(\nu) \geq -3$ and $F$ is not holomorphic
with respect to any of $I_1,\ I_2$ and $I_3$,
then $h^0(\eta) > 0$. But by Riemann-Roch,
\[
h^0(\eta) = e(\nu) - e(\xi) + 2-2g +
h^0((\num \oplus \xip) \otimes \kappa)
\]
where $\kappa$ is the canonical bundle of $\SI$.
Now $\xip \otimes \kappa$ is the trivial line bundle
and therefore, the last term is equal to $1+h^0(\num\otimes\kappa)$.
The result follows on noting that
$e(\xi) = 2g-2$ and that, for each $j \in \{1,2,3\},\
\pnm(I_jF_*(\rd/\rd z)) \otimes dz$ is a holomorphic section of
$\num \otimes \kappa$ which is nontrivial if
$F$ is not holomorphic with respect to $I_j$.
\end{pf}


\begin{thebibliography}{XXX}
\bibitem[AN]{an} Arezzo, C. and La Nave, G.,
Minimal two spheres and K\"{a}hler-Einstein metrics on Fano manifolds,
preprint available from
www.math.unipr.it/~rivista/membri/AREZZO/Arezzo.html
\bibitem[B-P-V]{bpv} Barth, W., Peters, C. and Van de Ven, A. ,
Compact Complex Surfaces, Springer-Verlag, Berlin 1984.
\bibitem[D]{d} Donaldson, S. K., Moment maps and diffeomorphisms,
Surv. Differ. Geom., VII, (2000), Int. Press, Somerville, MA, 107--127.
\bibitem[G-H]{gh} Griffiths, P. and Harris, J.,
Principles of Algebraic Geometry, John Wiley and Sons, New York 1978.
\bibitem[H]{h} Helgason, S.,
Differential Geometry, Lie Groups, and Symmetric Spaces,
Academic Press, 1978.
\bibitem[K]{ko} Kobayashi, R.,
Moduli of Einstein Metrics on a K3 surface and Degeneration of Type I,
Adv. Studies in Pure Math. {\bf 18-II} (1990) 257-311.
\bibitem[KT]{k} Kobayashi, R. and Todorov, N.,
Polarized period map for generalized K3 surfaces and
the moduli of Einstein metrics, Tohoku Math. J. {\bf 39} (1987) 341-363.
\bibitem[L]{l} Lawlor, G., The angle criterion,
Invent. Math. {\bf 95} (1989) 437-446.
\bibitem[M]{m} Micallef, M.,
Stable minimal surfaces in Euclidean space, J. Differential Geom.
{\bf 19} (1984), no. 1, 57-84.
\bibitem[Mo]{mo} Morgan, F., On the singular structure of
two-dimensional area minimizing surfaces in $\R^n$,
Math. Ann. {\bf 261} (1982) 101-110.
\bibitem[MW]{mw} Micallef, M. and Wolfson, J.,
The Second Variation of Area of Minimal Surfaces in Four-Manifolds,
Math. Ann. {\bf 295} (1993) 245-267.
\bibitem[SW]{sw} Schoen, R. and Wolfson, J.,
Minimizing areaamong lagrangian surfaces: the mapping problem,
J. Diff. Geom. {\bf 58} (2001) 1-86.
\bibitem[W]{w}  Wolfson, J., Lagrangian homology classes without regular 
minimizers, preprint.

\end{thebibliography}
\end{document}